\documentclass[12pt]{amsart}

\def\@typesizes{%
       \or{5}{6.5}\or{6}{7.5}\or{7}{8.5}\or{8}{11}\or{9}{12}%
       \or{10}{13}
       \or{\@xipt}{14}\or{\@xiipt}{15}\or{\@xivpt}{18}%
       \or{\@xviipt}{20}\or{\@xxpt}{24}}

\usepackage{amsthm}
\usepackage{amsmath}
\usepackage{amssymb}
\usepackage{graphicx}
\usepackage{enumerate}
\usepackage{color}

\allowdisplaybreaks
\numberwithin{equation}{section}
\numberwithin{figure}{section}

\theoremstyle{plain}
\newtheorem{theorem}{ Theorem}[section]
\newtheorem{proposition}[theorem]{ Proposition}
\newtheorem{lemma}[theorem]{ Lemma}
\newtheorem{corollary}[theorem]{ Corollary}
\newtheorem{example}[theorem]{ Example}
\newtheorem{remark}[theorem]{ Remark}
\newtheorem{definition}[theorem]{ Definition}
\newtheorem{conjecture}{ Conjecture}

\def\BET{\begin{theorem}}
\def\ENT{\end{theorem}}
\def\BEP{\begin{proposition}}
\def\ENP{\end{proposition}}
\def\BEL{\begin{lemma}}
\def\ENL{\end{lemma}}
\def\BEC{\begin{corollary}}
\def\ENC{\end{corollary}}
\def\BEE{\begin{example} \rm}
\def\ENE{\end{example}}
\def\BER{\begin{remark} \rm}
\def\ENR{\end{remark}}
\def\BED{\begin{definition} \rm}
\def\END{\end{definition}}
\def\BECJ{\begin{conjecture}}
\def\ENCJ{\end{conjecture}}

\def\bea{\begin{eqnarray}}
\def\eea{\end{eqnarray}}

\def\beas{\begin{eqnarray*}}
\def\eeas{\end{eqnarray*}}

\def\beq{\begin{equation}}
\def\eeq{\end{equation}}

\def\beal{\begin{align*}}

\def\eeal{ \end{align*} }

\def\roweq{\nonumber \\ &=& }
\def\rowleq{\nonumber \\  & \leq & }

\def\bfN{{\bf N}}

\def\bbC{{\mathbb C}}
\def\bbD{{\mathbb D}}

\def\bbN{{\mathbb N}}

\def\bbR{{\mathbb R}}

\def\cB{{\mathcal B}}

\def\ov{\overline}
\def\ef{\eqref}

\begin{document}

\title[Pointwise lower bounds]{Pointwise lower bounds in growth spaces with 
little $o$ conditions.}

\author{$\check{\mbox{Z}}$eljko $\check{\mbox{C}}$u$\check{\mbox{C}}$kovi\'c}
\address{Department of Mathematics, University of Toledo, Toledo, OH 43606, USA}

\author{Jari Taskinen}
\address{University of Helsinki, Department of Mathematics and Statistics,
P.O.Box 68, 00014 Helsinki, Finland}

\thanks{The second named author was partially supported by the V\"ais\"al\"a Foundation of the Finnish
Academy of Sciences and Letters and the EU-Horizon-MSCA project no. 101109510 LARGE BERGMAN.} 


\begin{abstract}
Pointwise lower bounds on the open unit disc
$\bbD$ for the sum of the moduli of two analytic functions $f$ and $g$ (or their derivatives)
are known in several cases, like $f,g$ belonging to the Bloch space $\cB$, $BMOA$ or the weighted Hardy space $H_\omega^\infty$. 
We modify the proofs of two important cases, proved by Ramey-Ullrich and Abakumov-Doubtsov, for functions
with little $o$ growth conditions. 
\end{abstract}

\maketitle

\section{Introduction}
\label{sec1}

The motivation for this paper comes from the important paper by Ramey and Ullrich \cite{RU}, where the authors proved the result about the growth of derivatives of functions in the Bloch space on the unit disk $\mathbb{D}.$  Subsequently several authors continued this line of investigation, in particular, we mention the work by Lou \cite{Lou}  studying 
the $\alpha$-Bloch spaces, Wulan and Zhu \cite{WZ} studying an analogue of the Ramey-Ullrich result for the unit ball in $\bbC^n$ and Abakumov and Doubtsov \cite{AD} who studied the weighted $H^{\infty}$ spaces; see also \cite{CR} and \cite{BLT} for
higher dimensional generalizations. The goal of our paper is to obtain complementary lower bounds for growth spaces with little $o$ condition.  In particular, we are interested in the little Bloch space, $VMOA$ and weighted $H^{\infty}$ spaces with $o$ growth conditions.  

We first introduce the relevant definitions. A function $g$ is called a Bloch function if it is analytic on $\mathbb{D}$ and it satisfies
$$ \sup_{z\in\mathbb{D} } (1 - |z|^2) |g'(z)| < \infty.$$
The space of Bloch functions is denoted by $\cB$. Similarly, the little Bloch space $\cB_0$ consists of functions $g$ analytic on $\mathbb{D}$ such that

$$ \lim_{|z|\rightarrow 1} (1 - |z|^2) |g'(z)| = 0.$$
Ramey and Ullrich proved the existence of two functions $f, g \in \cB$ such that
\bea
|f'(z)| + |g'(z)| \geq \frac{1}{1 - |z|}  \label{1.1}
\eea
for all $z\in \mathbb{D}.$
It should be noted that just one analytic function cannot satisfy the inequality above.  To prove their result, Ramey and Ullrich used functions that have a lacunary power series expansion in the form
$$
f(z) = \sum_{k=0}^{\infty} a_k z^{q^k}
$$
where $q$ is a positive integer.  It is also known that the function $f$ with such lacunary series belongs to the Bloch space $\cB$ if and only if $(a_k)_{k=0}^\infty \in \ell^{\infty}$. Similarly, $f\in \cB_0$ if and only if $(a_k)_{k=0}^\infty \in c_0$ (see \cite{Gir} or \cite{Po}).

Our Theorem \ref{th2.1} will give a growth condition for $\cB_0$ that complements the Ramey-Ullrich result.  In Proposition \ref{prop1.2} we give an analogous growth result for $VMOA$ functions.

Given  a radial weight function $\omega : \bbD \to (0,\infty)$, a generalized formulation of 
the question related to \ef{1.1} consists of finding  two analytic functions $f,g$ belonging to the weighted 
Hardy space $H^\infty_\omega$ and satisfying 
a pointwise lower bound 
\bea
|f(z)| + |g(z)| \geq \frac1{\omega (z)} .  \label{1.2}
\eea
A solution to  this problem  was found by Abakumov and Doubtsov in \cite{AD}, where they proved that
such functions can be found, if and only if the function $1/\omega$ is log-convex
on the interval $[0,1)$. Here, we will modify the proof of \cite{AD} so as to apply to functions
$f$ and $g$ with the corresponding small $o$ growth conditions. In that case one cannot expect to get
the same lower bound as in \eqref{1.2}, but we will show in Theorem \ref{th1.1} that if the right-hand side is 
multiplied by  an arbitrary continuous, radial  function  $W: \bbD \to (0,\infty)$ with $\lim_{|z|\to 1^-} W(z)=0$, then the functions $f$ and $g$ can be found. We also observe that this is in a sense the best possible lower bound for
$f,g$ with small $o$ growth conditions.

Here, by a radial weight $\omega$ we mean a continuous function on $\bbD$ 
with positive values, such that $\omega (z) = \omega(|z|) < \omega(|w|) =\omega(w)$ for all $z,w \in \bbD$ with $|z| > |w|$ and, moreover,
$\lim_{|z| \to 1^-}  \omega(z) = 0$. Then, the weighted Hardy space, or the growth space, 
$H_\omega^\infty$ is defined by 
\bea
H_\omega^\infty = \big\{  f: \bbD \to \bbC \mbox{ analytic} \, : \, \Vert f \Vert_\omega := 
\sup_{z \in \bbD} \omega(z) |f(z)| < \infty \big\},   \label{1.4}
\eea
and the corresponding space with $o$ growth condition is 
\bea
H_{\omega,0}^\infty = \big\{  f \in H_\omega^\infty \, : \, \lim_{|z| \to 1^-} \omega(z) |f(z)| =0  \big\}. 
\label{1.6}
\eea
As well known, both $H_\omega^\infty$ and $H_{\omega,0}^\infty$ are Banach spaces. Note that 
$f \in H_\omega^\infty$, if and only if, for some  constant $C_f > 0$ 
$$
|f(z)| \leq \frac{C_f}{\omega(z)} \ \ \ \  \forall\, z \in\bbD.
$$

We say that two positive, radial functions  $v$ and $w$ on $\bbD$ are equivalent, if there exists a constant $C>0$ 
such that $v(z) / C \leq w(z) \leq Cv(z)$ for all $z \in \bbD$. 
We say that $1/\omega$ is log-convex, if the function $- \log \omega(r)$ is equivalent with 
a convex function of $\log r$ for  $r \in (0,1)$. This holds 
in particular, if $- \log \omega (r)$ is a convex function of $r$, but there are other examples.
(See  Remark 6 in \cite{AD}. Note that in \cite{AD}, weights are defined as $1/\omega$, where $\omega$ 
is a weight in the sense of our paper.) 

Our modification of the main result in \cite{AD} reads as follows. 

\BET  \label{th1.1}
Assume that  $1/\omega$ is equivalent to a  log-convex function. 
Given any piecewise continuous, non-increasing function  $W : [0,1) \to (0, \infty)$ with the property
$\lim_{r \to 1^-} W(r) = 0$, there exist functions $f,g \in H_{\omega,0}^\infty$ such that
\bea
|f(z)| + |g(z)| \geq \frac{W(|z|)}{\omega(|z|)} \ \ \ \ \mbox{for all} \ z \in \bbD.   \label{3.6}
\eea
\ENT

Note that if $f,g \in H_{\omega,0}^\infty$, then the function $\widetilde W(z) : \bbD \to (0,\infty)$
\beas
\widetilde W(r) = \sup_{|w| \geq r} |f(w)|\omega(w) + \sup_{|w| \geq r} |g(w)|\omega(w) 
\eeas
is continuous, non-increasing and has the properties $\lim_{r \to 1^-} \widetilde W(r) = 0$
and 
$$
|f(z)| + |g(z)| \leq \widetilde W(|z|) / \omega(z) \ \ \ \forall \ z \in \bbD,
$$ 
by \ef{1.4}, \ef{1.6}. 
Hence, given two functions $f,g \in H_{\omega,0}^\infty$, there always exists a $W$ such that the reverse
inequality \ef{3.6} holds. In this sense the result of Theorem \ref{th1.1} gives the largest  possible lower 
bound for functions belonging to $H_{\omega,0}^\infty$.

\section{On the Ramey-Ullrich construction}
\label{sec2}

The result in  Proposition 5.4. of \cite{RU} is of course contained in Theorem 1.2. of \cite{AD},
nevertheless, we present in Theorem \ref{th2.1} a separate  formulation of our result in that case, too. 
The reason is that the proof is substantially simpler than that of Theorem \ref{th1.1}
and moreover, of its own interest, in view of the many applications of Proposition 5.4. of \cite{RU} appearing
in the literature. 

Our result 
reads as follows. 

\BET  \label{th2.1}
Given any continuous, non-increasing function  $W : [0,1) \to (0, 1 )$ with the property
$\lim_{r \to 1^-} W(r) = 0$, there exist functions $f,g \in \cB_0$ such that
\bea
|f'(z)| + |g'(z)| \geq \frac{W|z|)}{1 - |z|}  \ \ \ \mbox{for all} \ z \in \bbD.     \label{2.1}
\eea
\ENT

Proof. The  proof of Proposition 5.4. of \cite{RU} requires a number of
modifications. We assume  that the function $W$ is given and, for simplicity, write $q =100$, although
any large enough $q$ could be used in the same way as in the citation.

We define $f(z) = \sum_{j=1}^\infty a_j z^{q^j}$, where 
the numbers $a_j$ are defined inductively by
\bea
a_1 = 1 \ \ \ \mbox{and} \ a_j = \max \Big\{  W(1- q^{-j}), \frac{j-1}{j}  a_{j-1} \Big\}  
\ \forall \, j \geq 2. \label{2.4}
\eea
We first claim that 

\smallskip

\noindent $(i)$ the sequence $(a_j)_{j=1}^\infty$ is non-increasing, 

\smallskip

\noindent $(ii)$ $\lim_{j\to \infty } a_j=0$, and

\smallskip

\noindent $(iii)$ there holds $a_j \geq 1/j$ for all $j \in \bbN$.

\smallskip
 
To see $(i)$, we note that if $j \geq 2$ is arbitrary, then, 
$a_{j-1} \geq W(1 - q^{-(j-1)}) \geq W(1 - q^{-j}) $,  by observing the definition \ef{2.4} 
for the coefficient $a_{j-1}$ and the fact the $W$ is non-increasing. We obtain that 
the right hand side of \ef{2.4} is at most $a_{j-1}$. Thus, $a_j \leq a_{j-1}$, i.e., $(i)$ holds.  

Property $(iii)$ follows from the relations $a_1 =1$ and $a_j \geq (j-1)j^{-1} a_{j-1}$ by induction.  

 
In order to prove $(ii)$, we
denote by $\bfN \subset \bbN$ the set of indices
such that $a_j =  W(1- q^{-j})$, see \ef{2.4}. If $\bfN$ is finite, there is some number $J \in \bbN $ such that the relation 
$a_j = \frac{j-1}{j} a_{j-1}$ holds for all $j \geq J$ . By a repeated use of this relation
we obtain for all $\ell > J$
\bea
a_\ell = \Big( \prod_{j=J +1}^{\ell} \frac{j -1 }{j } \Big) a_J = \frac{J}{\ell} a_J \to 0 \ \ \ \mbox{as} \ \ell \to \infty.
\eea
If the set $\bfN$ is infinite and $\varepsilon > 0$ is given, we choose $J \in \bfN$ such
that $a_J =  W(1- q^{-J}) < \varepsilon$. (This is possible since $W(r) $ tends to 0 as $r \to 0$.)
By property $(i)$,  we have $a_j \leq a_J < \varepsilon$ for all $j > J$, which proves the claim.

The rest of the proof follows that of the citation, with some more  modifications. 
If $k \in \bbN$ is arbitrary, we first prove that 
\bea
|f'(z)| \geq \frac{CW(z)}{1- |z|} \ \ \ \mbox{for $z$ with} \ |z|
\in  I_k \subset [0,1), \label{2.10}
\eea
where we denote $ I_k = [1 - q^{-k} , 1 - q^{-(k +1/2)}]$, and to this end we use the triangle inequality to write
\bea
|z| \, |f'(z)| & \geq &   | a_k q^k z^{q^k}|  - \sum_{j=0}^{k-1} a_j q^j |z|^{q^j} - \sum_{j=k+1}^{\infty} a_j q^j |z|^{q^j}
\nonumber  \\
& =: &  I(z) - II(z) - III (z)  \label{2.12}
\eea
The estimate of the term I in \cite{RU} implies that 
\bea
I (z) \geq \frac{a_k q^k}3 
\ \ \ \mbox{for all} \ z \in I_k . \label{2.14}
\eea
Also, since the sequence $(a_j)_{j=1}^\infty$ is non-increasing, 
we have 
\bea
III(z) \leq  a_{k+1} \sum_{j=k+1}^{\infty} q^j |z|^{q^j} \leq  a_k \sum_{j=k+1}^{\infty} q^j |z|^{q^j}.  \label{2.16}
\eea
Now, in \cite{RU} it is shown that for $z \in I_k$ there holds
\bea
 \sum_{j=k+1}^{\infty} q^j |z|^{q^j} \leq  q^k \frac{q 2^{-q^{1/2}}}{1- q 2^{- q^{3/2} + q^{1/2} }}  . \label{2.18}
\eea
Our choice $q = 100$ implies 
\bea
q^{1/2} = 10  \Rightarrow q 2^{- q^{1/2}} = 100 \cdot 2^{-10} \leq 2^{-3}
\eea
and
\bea
q^{3/2} = 1000  \Rightarrow 1- q 2^{- q^{3/2} + q^{1/2}} \geq 1 - 100 \cdot 2^{-990} \geq 1 - \delta  
\eea
with $0 < \delta < 10^{-10}$. Hence, by \ef{2.16}, \ef{2.18}
\bea
III(z) \leq \frac{a_k q^k}{8(1- \delta)} \ \ \ \mbox{for all} \ z \in I_k .  \label{2.18a}
\eea

To treat the term $II(z)$ we denote by $\ell$ the  smallest integer not smaller than $k/2$ and write  
\bea
II(z) \leq \sum_{j=0}^{k-1} a_j q^j = \sum_{j=0}^{\ell-1} a_j q^j + \sum_{j=\ell}^{k-1} a_j q^j 
\ \ \ \mbox{for all} \ z \in I_k.  \label{2.19}
\eea
Since $0 < a_j < 1$ for all $j$, we can then estimate 
\bea
 \sum_{j=0}^{\ell-1} a_j q^j \leq  \sum_{j=0}^{\ell-1} q^j \leq \frac{q^{\ell }}{q-1}.   \label{2.20}
\eea
Since $\ell \leq k/2 +1$, $q =100$ and $a_k \geq 1/k$ (see $(iii)$), we obtain 
\bea
q^{\ell - k } < q^{- k/2 + 1 } < \frac{1}{100 k} \leq \frac{a_k}{100},
\eea
and combining this with \ef{2.20} yields
\bea
 \sum_{j=0}^{\ell-1} a_j q^j \leq   \frac{q^k q^{\ell -k}}{q-1} < \frac{a_k q^k}{100 (q-1)}.   \label{2.22}
\eea
Finally, the second sum in \ef{2.19} is estimated by 
\bea
\sum_{j=\ell}^{k-1} a_j q^j \leq a_{\ell} \sum_{j = \ell}^{k-1}  q^j
\leq a_{\ell} \sum_{j=0}^{k-1}  q^j \leq a_\ell \frac{q^k}{q-1} . \label{2.23}
\eea
Here, the choice \ef{2.4} implies
\bea
a_k \geq \Big( \prod_{j = \ell +1}^k  \frac{j-1}j \Big) a_\ell = \frac{\ell}{k} a_\ell \geq \frac{a_\ell}2 
\eea
which in combination with \ef{2.22} gives us 
\bea
\sum_{j=\ell}^{k-1} a_j q^j \leq a_k q^k\frac{2 }{q-1} . \label{2.24}
\eea
Altogether, \ef{2.18a}, \ef{2.19}, \ef{2.22}, and \ef{2.24} imply for $z \in I_k$ the bound
\bea
II(z) + III(z)  \leq a_k q^k \Big( \frac{3}{q-1}  + \frac{1}{8(1- \delta)} \Big) \leq \frac{ a_k q^k }{7} .
\eea
By \ef{2.12}, \ef{2.14}, 
\bea
|f'(z) | > |z| \,|f'(z)|  \geq I(z) - (II(z) + III(z) ) \geq \frac{ a_k q^k }{8}   \  \mbox{for} \ z \in I_k.   \label{2.26}
\eea
Finally, since $W$ is non-increasing, the maximum value of $W(z)$ on the interval $I_k$ 
is $W(1- q^{-k}) \leq a_k$ (see  \ef{2.10}, \ef{2.4}). Also, the  maximum
of the function $1/(1- |z|)$ on the interval $I_k$ equals  $q^k$. Thus, we have for  $z \in I_k$
\bea
\frac{W(|z|)}{1- |z|}  \leq  a_k q^k 
\eea
and \ef{2.26} implies $|f'(z)| \geq W(|z|)/(8(1-|z|))$ on $I_k$. 

Note that above, $k \in \bbN$ was arbitrary. Now, as in \cite{RU} and analogously to the above proof, one can 
construct a function $g$ such that
\bea
|g'(z)| \geq \frac{W(|z|)}{8(1- |z|)}
\eea
for all $z$ with $\ |z| \in [1- q^{- (k +1/2)} , 1-q^{-k-1} ]$. This and a suitable multiplication of the functions
$f$ and $g$ by a positive constant yields \ef{2.1}.\ \ 
$\Box$

\bigskip

Recall that a function $f \in L^1(\partial \bbD)$ belongs to the space $BMO$, if
\bea
\Vert f \Vert_* &:=& \sup_I \frac{1}{|I|} \int\limits_I \big| f(e^{i \theta}) - m_I (f) \big| d \theta < \infty,
\ \ \mbox{where} \ \ 
\nonumber \\ 
m_I (f) &=&  \frac{1}{|I|} \int\limits_I | f(e^{i \theta})| d \theta .  \label{2.27p}
\eea
Here, $I $ denotes an arbitrary interval with positive length contained in $\partial \bbD$.
The space $VMO$ consists of functions $f \in BMO$ such that the expression 
$|I|^{-1} \int_I | f(e^{i \theta}) - m_I (f) | d \theta$ 
tends to 0 as $|I| \to 0$. 
The space $BMOA$ consists of such Poisson extensions of functions in $BMO$ that are analytic
in $\bbD$, and the space $VMOA$ is defined in the same way with the help of $VMO$ instead of $BMO$.  
Endowed with the norm $\Vert f \Vert_{BMOA} = |f(0)| + \Vert f \Vert_*$, the space $BMOA$ is complete, and
$VMOA$ is its closed subspace. See \cite{Gir} or \cite{Zh} for more details.

\BEP  \label{prop1.2}
There exist functions $f,g \in VMOA$ such that
\bea
|f'(z)| + |g'(z)| \geq \frac{1}{(1 - |z|) \,\big| \log(1- |z|) \big|}
 \ \ \ \mbox{for all} \ z \in \bbD.     \label{2.30}
\eea
\ENP

Proof. We consider the function $W(r) = 1 / \big| \log(1- r) \big|$ in the proof 
of the previous theorem. Then, for every $j \in \bbN$, 
\bea
W(1- q^j) = \frac1{|\log q^{-j}|  } = \frac{1}{j \log 100} < \frac1j .
\eea
Thus, if we define $a_j = 1/j$ for all $j$, we see that the relations \ef{2.4} hold for all $j$. 
As a consequence of the proof of Theorem \ref{th2.1} and the relation \ef{2.1}, 
the functions 
\bea
f(z) =  \sum_{j=1}^\infty j^{-1} z^{q^j}  \label{2.32}
\eea
and $g$, constructed following \cite{RU}, satisfy the lower bound \ef{2.30}. 

By Theorem 9.3 in \cite{Gir}, $f$ and $g$ belong to $VMOA$ (as well as $H^2$), since the 
Taylor coefficients of the lacunary series \ef{2.32} belong to $\ell^2$. \ \ $\Box$

\section{Proof of Theorem \ref{th1.1}}
\label{sec3}

Proof of Theorem \ref{th1.1}. 
We modify the proof of Theorem 1.2 of \cite{AD}. Let us assume that we are given a weight $\omega$ 
such that $1/\omega$ is equivalent to a log-convex function. Our aim is to apply the considerations of Sections 
2.3--2.5 of \cite{AD} to the weight $w = 1/\omega$. In Section 2.3 of the citation the authors 
construct monotonically increasing sequences  $(x_k)_{k=0}^\infty$ and 
$(\beta_k)_{k=1}^\infty$ such that $x_k < 0 $ and $\beta_k >0$ for all $k \in \bbN$ and
$\lim_{k \to \infty} x_k = 0$ and $\lim_{k \to \infty} \beta_k = \infty$, and denote 
\bea
t_k = \exp (x_k) \in (0,1) \ \ \ \ \mbox{for all $k$.}   \label{3.7}
\eea
They also define, for all $k \in \bbN$, certain parameters $a_k > 0$ and linear functions 
$\ell_k: (-\infty,0) \to \bbR$ 
\bea
\ell_k(x) = \log a_k + \beta_k x , \ \ \ \mbox{hence}, \ \ e^{\ell_k(x)} = a_k r^{\beta_k},  \label{3.8}
\eea
where $x \in (-\infty,0),$  $r = e^x  \in (0,1)$. Note that in the construction, the numbers
$x_0 < 0$ and $h >0$ can be chosen arbitrarily, but after that all other parameters are uniquely 
determined; see Section 2.3. and the beginning of Section 2.5 of \cite{AD}. 
Here, we fix $h \geq 8$ for the rest of the proof (in the reference, the results hold for all $h \geq 2$).
The number $x_0$ will be fixed later.

We do not need to present the details of the above mentioned construction, since it
is enough for our purposes to quote a couple of the results in \cite{AD}.  First,  the inequality (2.2) in Lemma 2.1. of \cite{AD} says that 
\bea
\ell_{k+1}(x)  \geq \ell_k(x) + h \ \ \ \mbox{for all} \ x \in [x_{k+1} ,0 ), \ k \in\bbN,
\label{3.10}
\eea
and the choice of the functions $\ell_k$ (see (ii) in the beginning of Section 2.3 and Fig.1 of \cite{AD})
implies that 
\bea
\ell_{k+1}(x)  \geq \ell_k(x)  \ \ \ \mbox{for all} \ x \in [x_{k} , x_{k+1} ), \ k \in\bbN,
\label{3.10a}
\eea
Second, Lemma 2.2 of \cite{AD}  states that  the following holds:

\bigskip

\noindent $(I)$ \ \ $a_k r^{\beta_k} \leq 1/ \omega(r), \ \ \ r \in [t_0,1), $

\smallskip

\noindent $(II)$  \ \ $e^{-h} /\omega(r) \leq a_k r^{\beta_k} , \ \ \ r \in [t_{k-1}, t_k], $

\smallskip

$$
\sum_{\stackrel{\scriptstyle m \geq 1, }{|m-k| \geq 2}
}  a_m r^{\beta_m}    \leq \frac12 a_k r^{\beta_k} , \ \ \ r \in [t_{k-1}, t_k]. \leqno(III)
$$

If necessary, we apply a simple scaling argument so that we can assume for the rest of the proof that $W(0) =1$ and thus $W(|z|) \leq 1$ for all $z \in \bbD$. 
We define $\nu_1=1$ and 
\bea
\nu_k = \max \Big\{  W(t_{k}), \frac{k-1}{k}  \nu_{k-1} \Big\}    \label{3.11}
\ \forall \, k \geq 2.
\eea
Then, in the same way as after \ef{2.4}, one shows that the sequence $(\nu_k)_{k=1}^\infty$
is non-increasing and there holds $\lim_{k \to \infty} \nu_k = 0$. Also, 
we have $\nu_k \geq 1/k$ for all $k$, and, moreover, if $m < k$, there holds, by a repeated use of \ef{3.11}, 
\bea
\nu_m \leq \Big( \prod_{j=m+1}^k  \frac{j}{j-1} \Big) \nu_k = \frac{k}{m} \nu_k . \label{3.11a}
\eea

Due to the properties of the numbers $\nu_k$, the function  
$\widehat{W} : (0,1) \to (0,1)$, 
\bea
\widehat{W} (r) = 1 \  \mbox{for} \ r \in [0,t_0),
\ \ \ \widehat{W} (r) = \nu_k  \ \mbox{for} \ r \in [t_{k}, t_{k+1}), \ k \in \bbN_0,
\label{3.11b}
\eea
is non-increasing, piecewise continuous and has the limit 0 as $r \to 1^-$.
Moreover, since $\widehat{W}$ is not smaller than $W$, it suffices to prove the assertion of the  
theorem for this new function, which we still denote by $W$. 

Our next task is to prove the following  

\bigskip

{\it Claim.} 
There holds the following modification of $(III)$: for arbitrary  $k \in \bbN$ we have 
$$
\sum_{\stackrel{\scriptstyle m \geq 1, }{|m-k| \geq 2}
} \nu_m a_m r^{\beta_m}    \leq \frac12 \nu_k a_k r^{\beta_k} , \ \ \ r \in [t_{k-1}, t_k],
\leqno(IV) 
$$

Now, the main difficulty is that if $m < k$, then the new coefficient $\nu_m$, on the left hand side, 
is larger than $\nu_k$ (thus, the sequence $(\nu_k)_{k=1}^\infty$ must not decrease too fast).
To prove the claim we first consider indices $m$ smaller than $k$ on the left-hand side of $(IV)$.
Indeed, \ef{3.10} implies $\ell_{j+1}(x) - \ell_j(x) \geq h$ for all  $j = m ,m +1, \ldots , k-2$
and $x \in [x_{k-1} , 0)$, hence, for these $x$,
\bea
& & \ell_m(x) + h(k-m-1) \leq \ell_m (x) + \sum_{j=m}^{k-2} \big( \ell_{j+1}(x) - \ell_j(x) \big)  
\roweq 
\ell_{k-1}(x) \leq \ell_k(x)  ,  \label{3.14}
\eea
see also \ef{3.10a}. 
We obtain, by \ef{3.8}, \ef{3.11a} and \ef{3.14}, for all $m \leq k-2$ and $r = e^x \in [t_{k-1},1) $ 
(see \ef{3.7}),
\bea
&&   \frac{m}{k} e^{h(k-m-1)}  \nu_m a_m r^{\beta_m} = \nu_m \frac{m}{k}
e^{ \ell_m(x)}  e^{h(k-m-1)} \leq \nu_k e^{ \ell_m(x) + h(k-m-1)}
\rowleq  
\nu_k e^{ \ell_k(x) }= \nu_k a_k r^{\beta_k}  . 
\eea
Hence,
\bea
& & \sum_{m = 1}^{k-2}  \nu_m a_m r^{\beta_m} \leq \nu_k a_k r^{\beta_k}
\sum_{m = 1}^{k-2} \frac{k}{m} e^{- h(k-m-1)} . 
\label{3.16}
\eea
Here, assuming $k \geq 3$, $h \geq 8$, 
\bea
& & \sum_{m = 1}^{k-2} \frac{k}{m} e^{- h(k-m-1)} 
\roweq
e^{- h(k-1)}\sum_{m = 1}^{[k/2]} \frac{k}{m} e^{hm} + e^{-h}
\sum_{m = [k/2]+1}^{k-2} \frac{k}{m} e^{- h(k - 2 - m)}
\rowleq
e^{- h(k-1)} k^2 e^{h k / 2} + \frac{k}{k/2}e^{-h} \sum_{m = 0}^{[k/2] - 2 }  e^{-hm}
\rowleq 
e^{- h(k/2 - 1)} k^2 +  2 e^{- h } \sum_{m=0}^\infty e^{- hm}.
\eea
The first term on the right is at most $ 
e^{- 4k +8 } k^2 $ (since $h \geq 8$), which is decreasing for $k \geq 3$ and thus
at most  $9 e^{-4} <  1/5$. The second term is bounded by
$ 3 e^{-h} < \frac1{20}. $  Hence,  \ef{3.16} yields
\bea
& & \sum_{m = 1}^{k-2}  \nu_m a_m r^{\beta_m} \leq \frac14 \nu_k a_k r^{\beta_k} \label{3.18}
\eea
On the other hand, for $m > k$  we can use the estimate just before (2.3) of \cite{AD}, namely,
\bea
a_m r^{\beta_m} \leq a_k r^{\beta_k} e^{-h (m-k-1)} , \ \ \ r \in [t_{k-1} ,t_k], 
\eea
and since $\nu_k > \nu_m$ for $m > k$, this yields 
\bea
\sum_{m = k+2}^\infty  \nu_m a_m r^{\beta_m}  < \nu_k 
\sum_{m = k+2}^\infty a_m r^{\beta_m} 
\leq \nu_k a_k r^{\beta_k} \sum_{m=1}^\infty e^{- mh}
\leq \frac14 \nu_k a_k r^{\beta_k}    \label{3.20}
\eea
for $r \in [t_{k-1}, t_k].$
Combining this with \ef{3.18} proves $(IV)$. \ \ $\boxtimes$ 

\bigskip

The rest of the proof follows the main reference with some necessary changes. Accordingly, we denote $n_k = [\beta_k] + 1 \geq \beta_k$ for $k \in\bbN$.
The arguments in the beginning of Section 2.5 in \cite{AD} include the choice of the
number $x_0 < 0$ such that $t_0 = e^{x_0} > 9/10$, which then implies  that $n_k < n_{k+1}$
and that 
\bea
1 < \frac{a_k r^{\beta_k} }{a_k r^{n_k} } \leq \frac{10}{9}
\eea
holds for all $r$ close enough to 1, that is, at least for $r \in (\frac{9}{10},1)$. Thus, 
properties $(I)$, $(II)$, $(IV)$  turn into

\bigskip

\noindent $(A)$ \ \ $a_k r^{n_k} \leq 1/ \omega(r), \ \ \ r \in [t_0,1), $

\smallskip

\noindent $(B)$  \ \ $\frac{9}{10} e^{-h} /\omega(r) \leq a_k r^{n_k} , \ \ \ r \in [t_{k-1}, t_k], $

\smallskip

$$
\sum_{\stackrel{\scriptstyle m \geq 1, }{|m-k| \geq 2}
} \nu_m a_m r^{n_m}    \leq \frac12\cdot \frac{10}{9} \nu_k a_k r^{n_k} , \ \ \ r \in [t_{k-1}, t_k], \leqno(C)
$$

We define
\bea
f(z) = \sum_{j=0}^\infty \nu_{2j+1} a_{2j +1} z^{n_{2j+1}} , \ \ \ 
g(z) = \sum_{j=1}^\infty \nu_{2j} a_{2j} z^{n_{2j}}  \label{3.22}
\eea
and claim first that both functions belong to $H_{\omega,0}^\infty$. Indeed 
for every $k \in \bbN$ and $z$ with $r = |z| \in [t_{k-1}, t_k]$ we obtain
\bea
& & |f(z)| + |g(z)| \leq \sum_{m=1}^\infty \nu_m a_m |z|^{n_m}
\roweq 
\sum_{m=k-1}^{k+1} \nu_m a_m r^{n_m} + 
\sum_{\stackrel{\scriptstyle m \geq 1, }{|m-k| \geq 2}} \nu_m a_m r^{n_m}. \label{3.24}
\eea
By $(A)$ and \ef{3.11a}, the first sum on the right-hand side can be estimated as 
\beas
 \sum_{m=k-1}^{k+1} \nu_m a_m r^{n_m} \leq \sum_{m=k-1}^{k+1}
\frac{\nu_m }{\omega(r)} \leq 3\frac{(k+1)\nu_{k+1} }{(k-1)\omega(r)}
\leq \frac{6 \nu_{k+1} }{\omega(r)}.
\eeas
Since $r  \in [t_{k-1}, t_k]$, we have $W(r) \geq \nu_{k-1} >  \nu_{k+1}$ by \ef{3.11b},
and the first sum in \ef{3.24} thus has the upper bound $6 W(r) /\omega(r)$.  
The second sum in \ef{3.24} is also bounded by $CW(r)/\omega(r)$, since $(C)$ and similar arguments yield
\beas
 \sum_{\stackrel{\scriptstyle m \geq 1, }{|m-k| \geq 2}} \nu_m a_m r^{n_m} \leq 
\frac12\cdot \frac{10}{9} \nu_k a_k r^{n_k}  \leq \frac{\nu_k }{ \omega(r)}
\leq \frac{W(r) }{\omega(r)}.
\eeas

Hence, we obtain 
\beas
|f(z)| + |g(z)| \leq  \frac{CW(|z|)}{\omega(z)}
\eeas
for all $z \in \bbD$, since $k$ was arbitrary. This proves the claim since $W(r) \to 0$ as $r \to 1^-$.

To prove the lower bound \ef{3.6} we consider first an odd $k \in \bbN$, $k \geq 3$, and $z \in \bbD $ with 
$|z| = r \in [t_{k-1}, t_k)$. Using $(B)$, $(C)$ and the bound
$
W(r)  = \nu_{k-1} \leq k(k-1)^{-1} \nu_k 
$  (see \ef{3.11b}, \ef{3.11a})
we estimate
\bea
& & e^{-h} \frac{W(r)}{\omega(r)} \leq \frac{10k}{9(k-1)} \nu_k a_k r^{n_k}
\leq 2 \nu_k a_k r^{n_k} 
\rowleq 
5 \bigg( \nu_k a_k r^{n_k} - \!\!\!\!\! \sum_{\stackrel{\scriptstyle m \geq 1, }{|m-k| \geq 2}} 
\!\!\!\! \nu_m a_m r^{n_m} \bigg)
= 
5 \bigg( \big|\nu_k a_k z^{n_k} \big|- \!\!\!\!\! \sum_{\stackrel{\scriptstyle m \geq 1, }{|m-k| \geq 2}} 
\!\!\!\! \nu_m a_m |z|^{n_m}  \bigg)
\rowleq
5  \bigg| \nu_k a_k z^{n_k} + \!\!\!\!\! \sum_{\stackrel{\scriptstyle m \geq 1, }{|m-k| \geq 2}} 
\!\!\!\! \nu_m a_m z^{n_m} \bigg| = 5 |f(z)| .
\eea

If $k \in \bbN$ is even and $r \in [t_{k-1}, t_k]$, the same argument yields
the lower bound $5 g(r) \geq e^{-h} W(r) \omega(r)^{-1}$. By a simple redefinition of 
$f$ and $g$, this proves \ef{3.6} for all $z$ belonging to the annulus $\bbD \setminus \ov{B(0,t_2)}$. 

Condition \ef{3.6} could still be violated  in case the functions $f$ and $g$  have 
common zeros in $B(0,t_2)$.  The common zero of $f$ and $g$ at the origin is treated 
by dividing $f$ by $z^{n_1}$ and redefining $f$ accordingly. 
Then, the lower bound \ef{3.6} can be obtained in the entire disc $\bbD$ by choosing
a suitable $\theta \in (0,2\pi)$ so that changing the definition of $g$ by the 
rotation $g(z) \mapsto g(e^{i \theta}z)$ 
removes the possible remaining, finitely many common zeros of $f$ and $g$. For details, 
see the end of Section 2.5 of \cite{AD}.  \ \ $\Box$

\end{document}